
\baselineskip=14pt
\parskip=10pt

\magnification=\magstephalf

\def\P{{\cal P}}

\def\1{{\overline{1}}}
\def\2{{\overline{2}}}
\parindent=0pt
\overfullrule=0in

\def\frac#1#2{{#1 \over #2}}

\bf
\centerline
{
The C-finite Ansatz meets the Holonomic Ansatz
}

\rm
\bigskip
\centerline
{\it By Shalosh B. EKHAD and Doron ZEILBERGER}
\bigskip
\qquad 

{\bf VERY IMPORTANT} 

As in all our joint papers, the main point is not the article, but the accompanying Maple package, {\tt CfiniteIntegral.txt}, that
may be downloaded, free of charge, from the web-page of this article

{\tt http://www.math.rutgers.edu/\~{}zeilberg/mamarim/mamarimhtml/cfiniteI.html} \quad ,

where the readers can also find sample input and output files, that they are welcome to extend using their
own computers.

{\bf Preface}

In a recent article, [Kim], symbolic summation, and quite a bit of human pre-processing, is  used to evaluate
certain integrals involving Chebyshev polynomials.

First, let's remark that the Chebyshev polynomials, like all {\it classical orthogonal polynomials},
belong to the {\bf Holonomic Ansatz} ([Z1], beautifully, and very efficiently, implemented in [Kou]),
and as such, {\it inter alia}, {\bf every} identity in [Kim] (and in many other articles that are still
published today) are {\it automatically provable}, and their {\it epistemological status} is the same as
identities like $(a^3-b^3)=(a-b)(a^2+2ab+b^2)$ or $134 \cdot 431=57754$. Could you imagine a paper
published today (or even two thousands years ago), entitled ``A new proof of the identity  $134 \cdot 431=57754$''?,

{\bf Introduction}

But the Chebyshev polynomials are not `just' {\it holonomic}, they belong to the more restricted class
of $C$-finite polynomial sequences, and hence belong to the $C$-finite ansatz ([Z2]), and as such
have nice closure properties. It turns out that one can {\it interface} the $C$-finite ansatz and
the {\it holonomic} ansatz, and borrow from the latter the powerful, and not-as-well-known-as-it-should-be

``{\bf the (continuous) Almkvist-Zeilberger algorithm}'',

described in [AZ], and implemented in the Maple package

{\tt http://www.math.rutgers.edu/\~{}zeilberg/tokhniot/EKHAD.txt} \quad .

We combined all the procedures from the above Maple package and the Maple package that accompanies [Z2], and
created a new, self-contained, Maple package,

{\tt http:/www.math.rutgers.edu/\~{}zeilberg/tokhniot/CfiniteIntegral.txt} \quad ,

that can {\bf automatically} prove {\it every} identity in [Kim], and many, {\bf far deeper}, identities.

We should note that [Kim] also uses computer algebra methods, but those involving {\it summation}, and spends
quite a lot of {\it human effort} to  go from integration to summation. With the
Almkvist-Zeilberger algorithm, one can proceed directly, as follows.

{\bf Using the Continuous Almkvist-Zeilberger algorithm in order to Evaluate (Symbolically!)  Integrals of Powers of [in particular] Chebyshev Polynomials}

Suppose that you want to study a sequence of the form
$$
a(n):= \int_{\alpha}^{\beta} P_n(x) \, K(x) \,dx \quad,
$$
for some ``nice'' {\it kernel}, $K(x)$, and a sequence $P_n(x)$ of $C$-finite polynomials, i.e. given by a recurrence
$$
P_n(x)= \sum_{i=1}^{L} p_i(x) P_{n-i}(x) \quad ,
$$
for some positive integer, $L$, (the {\it order}) and polynomials $p_i(x)$ ($i=1, \dots ,L$), subject to
{\it initial conditions}
$$
P_0(x)=q_0(x), \dots , P_{L-1}(x)=q_{L-1}(x) \quad ,
$$
for some polynomials $q_0(x), \dots, q_{L-1}(x)$.

What's nice about the C-finite ansatz is that once you know that $\{P_n(x)\}$ is $C$-finite, the same is true
for  $\{P_n(x)^r\}$ for any positive integer $r$,
and also the sequence $P_n(x)P_n^{*}(x)$ (where for any polynomial $a(x)$,
$a^{*}(x)=x^d a(1/x)$ (where $d$ is the degree of $a(x)$), is its ``reverse''), and many other related sequences,
and one can fully automatically (and very fast) find C-finite representations for them (see [Z2], and  [KP] (a true {\bf masterpiece}!)).

Once we have such a $C$-finite polynomial sequence, the {\it ordinary} generating function
$$
R(x,t):=\sum_{n=0}^{\infty} P_n(x) t^n \quad,
$$
is a certain {\it rational function}, $R(x,t)$, of the variables $x$ and $t$, that can be easily, and quickly, found automatically. 
Hence the (ordinary) generating function of the sequence $a(n)$, let's call it $f(t)$, 
$$
f(t):=\sum_{n=0}^{\infty} a(n) \, t^n \quad,
$$
can be expressed as
$$
f(t)= \int_\alpha^\beta \, R(x,t) K(x) \, dx \quad .
$$
The continuous Almkvist-Zeilberger algorithm  produces a  {\it linear differential operator} with {\it polynomial} coefficients,
$\P(t,\frac{d}{dt})$, and a {\it certificate} (a rational function times the integrand), let's call it $C(x,t)$, such that
$$
\P(t,\frac{d}{dt})  \left [ \, R(x,t) K(x) \, \right ] \, = \, \frac{d}{dx} C(x,t) \quad .
$$
Integrating with respect to $x$,  from $x=\alpha$ to $x=\beta$, we get
$$
 \P(t,\frac{d}{dt}) \left [ \, \int_\alpha^\beta \,  R(x,t) K(x)\, dx \, \right ]  = \int_\alpha^\beta \, \left ( \frac{d}{dx} C(x,t) \right ) \quad .
$$
Hence, by the {\it Fundamental Theorem of Calculus}, $f(t)$ satisfies an {\it inhomogeneous} ordinary differential
equation with polynomial coefficients
$$
 \P(t,\frac{d}{dt}) f(t)= C(\beta,t) - C(\alpha,t) \quad .
\eqno(DiffEq)
$$
It is readily seen that the right-hand-side is a rational function in $t$. The differential equation $(DiffEq)$ immediately
implies a {\it linear (inhomogeneous) recurrence equation with polynomial coefficients} for the actual coefficients,
namely, our original $a(n)$, from which, using
standard methods, one can get a {\it homogeneous} linear recurrence equation with polynomial coefficients,
that together with the {\it initial conditions}, that can be easily computed, constitutes a full
{\bf description} of the sequence $a(n)$, that enables us to compute the sequence to as-many-as-desired terms.

In fact, since we have the {\it theoretical  guarantee} that such a description {\bf exists} we can even
skip the above steps, and resort to {\it pure guessing}!

As we have mentioned at the beginning, the above-mentioned web-page

{\tt http://www.math.rutgers.edu/\~{}zeilberg/mamarim/mamarimhtml/cfiniteI.html} \quad ,

contains numerous sample input and output files (including fully automatic proofs of the results of [Kim]),
that readers can extend to their heart's content.

{\bf References}

[AZ]  Gert Almkvist and Doron Zeilberger,
{\it The Method of Differentiating Under The Integral Sign}, J. Symbolic Computation {\bf 10}(1990), 571-591. Available from  \hfill\break
{\tt http://www.math.rutgers.edu/\~{}zeilberg/mamarim/mamarimhtml/duis.html} \quad .

[KP] Manuel Kauers  and Peter Paule, {\it ``The Concrete Tetrahedron''}, Springer, 2011.

[Kim] Seon-Hong Kim, {\it On some integrals involving Chebyshev polynomials},
Ramanujan Journal {\bf 38} (2015), 629-639.

[Kou] Christoph Koutschan, {\it Holonomic functions in Mathematica}, ACM Communications in Computer Algebra {\bf 47}(2013), 179-182. Available from  \hfill\break
{\tt http://www.koutschan.de/publ/Koutschan13b/holofunc.pdf} \quad .

[Z1] Doron Zeilberger,  {\it A Holonomic Systems Approach To Special Functions}, 
J. Computational and Applied Math {\bf 32} (1990), 321-368. Available from \hfill\break
{\tt http://www.math.rutgers.edu/\~{}zeilberg/mamarim/mamarimhtml/holonomic.html} \quad .

[Z2] Doron Zeilberger, {\it The C-finite Ansatz},  Ramanujan Journal {\bf 31} (2013), 23-32. Available from  \hfill\break
{\tt http://www.math.rutgers.edu/\~{}zeilberg/mamarim/mamarimhtml/cfinite.html} \quad .

\bigskip
\bigskip
\hrule
\bigskip
Doron Zeilberger, Department of Mathematics, Rutgers University (New Brunswick), Hill Center-Busch Campus, 110 Frelinghuysen
Rd., Piscataway, NJ 08854-8019, USA. \hfill \break
zeilberg at math dot rutgers dot edu \quad ;  \quad {\tt http://www.math.rutgers.edu/\~{}zeilberg/} \quad .
\bigskip
\hrule
\bigskip
Shalosh B. Ekhad, c/o D. Zeilberger, Department of Mathematics, Rutgers University (New Brunswick), Hill Center-Busch Campus, 110 Frelinghuysen
Rd., Piscataway, NJ 08854-8019, USA.
\bigskip
\hrule

\bigskip
Exclusively published in The Personal Journal of Shalosh B. Ekhad and Doron Zeilberger  \hfill \break
({ \tt http://www.math.rutgers.edu/\~{}zeilberg/pj.html})
and {\tt arxiv.org} \quad . 
\bigskip
\hrule
\bigskip
{\bf Dec. 21, 2015}

\end